\documentclass[11pt]{amsart}

\usepackage{amsfonts}
\usepackage{amssymb}
\newcommand\nc{\newcommand}
\input xy
\xyoption{all}

\newtheorem{theorem}{Theorem}[section]
\newtheorem{prop}[theorem]{Proposition}
\newtheorem{prblm}[theorem]{Problem}
\newtheorem{defin}[theorem]{Definition}
\newtheorem{caution}[theorem]{Caution}
\newtheorem{jtheorem}[theorem]{Theorem}

\newtheorem{remark}[theorem]{Remark}
\newtheorem{restroom}[theorem]{Restroom}
\newtheorem{lemma}[theorem]{Lemma}
\newtheorem{construction}[theorem]{Construction}
\newtheorem{corollary}[theorem]{Corollary}
\newtheorem{example}[theorem]{Example}
\newtheorem{conclusion}[theorem]{Conclusion}
\newtheorem{triviality}[theorem]{Triviality}
\newtheorem{proto}[theorem]{Prototype Quasifibration}
\newtheorem{cauex}[theorem]{Cautionary Example}
\newtheorem{propositiondef}[theorem]{Proposition-Definition}
\newtheorem{subth}{Nuisance}[theorem]
\newtheorem{ssubth}{ }[subth]
\newtheorem{conjecture}[theorem]{Conjecture}
\newtheorem{question}[theorem]{Question}

\nc\ques[1]{\begin{question}
\label{#1}
\begin{em}}
\nc\tri[1]{\begin{triviality}
\label{#1}}
\nc\rest[1]{\begin{restroom}
\label{#1}}
\nc\conj[1]{\begin{conjecture}
\label{#1}}
\nc\prodef[1]{\begin{propositiondef}
\label{#1}}
\nc\prt[1]{\begin{proto}
\label{#1}}
\nc\lem[1]{\begin{lemma}
\label{#1}}
\nc\pro[1]{\begin{prop}
\label{#1}}
\nc\thm[1]{\begin{theorem}
\label{#1}}
\nc\teiri[1]{\begin{jtheorem}
\label{#1}}
\nc\cor[1]{\begin{corollary}
\label{#1}}
\nc\dfn[1]{\begin{defin}
\label{#1}}
\nc\sthm[1]{\begin{subth}
\label{#1}}
\nc\exm[1]{\begin{example}
\label{#1}
\begin{em}}
\nc\plm[1]{\begin{prblm}
\label{#1}
\begin{em}}
\nc\rmk[1]{\begin{remark}
\label{#1}
\begin{em}}
\nc\ntn[1]{\begin{notation}
\label{#1}
\begin{em}}
\nc\cau[1]{\begin{caution}
\label{#1}
\begin{em}}
\nc\imn[1]{\begin{importnota}
\label{#1}
\begin{em}}
\nc\cax[1]{\begin{cauex}
\label{#1}
\begin{em}}
\nc\con[1]{\begin{construction}
\label{#1}
\begin{em}}
\nc\ssthm[1]{\begin{ssubth}
\label{#1}
\begin{em}}
\nc\cnc[1]{\begin{conclusion}
\label{#1}
\begin{em}}

\nc\eques{\end{em}
\end{question}}
\nc\elem{\end{lemma}}
\nc\erest{\end{restroom}}
\nc\econj{\end{conjecture}}
\nc\eprodef{\end{propositiondef}}
\nc\eprt{\end{proto}}
\nc\ethm{\end{theorem}}
\nc\eteiri{\end{jtheorem}}
\nc\ecor{\end{corollary}}
\nc\edfn{\end{defin}}
\nc\esthm{\end{subth}}
\nc\epro{\end{prop}}
\nc\etri{\end{triviality}}
\nc\eexm{\end{em}
\end{example}}
\nc\ermk{\end{em}
\end{remark}}
\nc\eplm{\end{em}
\end{prblm}}
\nc\ecau{\end{em}
\end{caution}}
\nc\ecax{\end{em}
\end{cauex}}
\nc\eimn{\end{em}
\end{importnota}}
\nc\entn{\end{em}
\end{notation}}
\nc\econ{\end{em}
\end{construction}}
\nc\ecnc{\end{em}
\end{conclusion}}
\nc\essthm{\end{em}
\end{ssubth}}

\nc\GS{{\mathfrak S}}
\nc\pp{{\mathbb P}}
\nc\qq{{\mathbb Q}}
\nc\Q{{\mathbb Q}}
\nc\oo{{\mathcal O}}
\nc\sss{{\mathbb T}}
\nc\zz{{\mathbb Z}}
\nc\Z{{\mathbb Z}}
\nc\rr{{\mathbb R}}
\nc\R{{\mathbb R}}
\nc\cc{{\mathbb C}}
\nc\C{{\mathbb C}}
\nc\ff{{\mathbb F}}
\nc\ee{{\mathbb E}}
\nc\AAA{{\mathbb A}}
\nc\disjoint{{\amalg}}
\nc\id{{\rm id}}
\nc\vv{{\mathbb V}}
\nc\Ker{{\rm Ker }}
\nc\sgn{{\rm sgn}}
\nc\CCC{{\mathcal C}}
\nc\unit{{\hbox{\bf 1}}}
\nc\Poin{{\mathcal P}}
\nc\FFF{{\mathcal F}}
\nc\pic{{\mathcal Pic}}
\nc\Pic{{\rm Pic}}
\nc\alb{{\mathcal Alb}}
\nc\AVar{{\mathcal Var}/\kappa}
\nc\ChM{{\mathcal ChM}}
\nc\alphabar{{\overline{\alpha}}}
\nc\gm{{{\mathbb G}_m}}
\nc\Hodge{{\mathcal HS}}
\nc\Vect{{\mathcal Vec^{\pm}}}

\author{E. Javier Elizondo}
\thanks{The first author is partially supported by grants DGAPA 107012 and CONACYT 101519}
\author{Shun-ichi Kimura}
\thanks{The second author is partially supported by the program JSPS-CONACYT}

\title[Motivic Chow Series]{Rationality of Motivic Chow Series modulo
  $\AAA^1$-homotopy} 


\begin{document}

\begin{abstract}
Consider the formal power series $\sum [C_{p, \alpha}(X)]t^{\alpha}$ (called
Motivic Chow Series), where $C_p(X)=\disjoint C_{p, \alpha}(X)$ is the Chow
variety of $X$ parametrizing the $p$-dimensional effective cycles on $X$
with $C_{p, \alpha}(X)$ its connected components, and $[C_{p, \alpha}(X)]$
its class in $K(\ChM)_{\AAA^1}$, the $K$-ring of Chow motives modulo
$\AAA^1$ homotopy.  Using Picard product formula
 and Torus action, we will show that the Motivic Chow Series
is rational in many cases.
\end{abstract}
\maketitle
\tableofcontents

\setcounter{section}{-1}
\section{Introduction}

Let $X$ be an algebraic variety over a finite field $\ff_q$.  Andr\'e Weil
\cite{W} conjectured that the formal power series  
$\sum_d |Sym^d(X)(\ff_q)|t^d$
is a rational function, in this paper we call this series as the Weil
zeta series. This is part of what it is known as 
the Weil conjecture and was proved by Dwork  in 1960, see \cite{D}.
Weil also observed that 
if a suitable cohomology theory exists, axiomatized as  
Weil cohomology, then the rationality of the Weil zeta series
follows.  Furthermore, the degree of the denominator is the 
dimension of the odd part of the cohomology ring, and the degree of
the numerator is the dimension of the even part.   

In 2000, Kapranov \cite{Ka} proved that if $X$ is a smooth projective
curve, then the series  $\sum [Sym^dX]t^d$, that we call
Kapranov motivic series, is a rational function in 
$K'(\AVar)[[t]]$, where $K'(\AVar)$ is the $K$-ring of  algebraic
varieties over a fixed field $\kappa$, together with the relation
$[X]=[C]+[X-C]$ where $C\subset X$ is a closed subscheme.  The
rationality of the Kapranov motivic zeta series
implies the rationality of the Weil motivic zeta series when the base
field is a finite 
field.  Moreover, it makes sense to ask if the 
Kapranov motivic zeta is rational for a general variety, even when the
base field $\kappa$ is an infinite field.  

In 2004, Larsen and Lunts  \cite{LL} proved that when $X$ is a
surface, the Kapranov motivic zeta series $\sum [Sym^dX]t^d$ is rational in
$K'(\AVar)[[t]]$ if and only if $X$ is a ruled surface. This tells us
that we can 
expect the rationality of the Kapranov motivic zeta in $K'(\AVar)$
only for special varieties. 

Then in 2005 Y. Andr\'e  \cite{A} observed that if the motive of $X$ is finite
dimensional, then the formal power series $\sum [Sym^dX]t^d$ in
the  ring of formal power series with coefficients in the $K$-ring of Chow motives,
denoted by $K(\ChM)[[t]]$, is rational, as before we call this series
the  Andr\'e motivic zeta series.  The Chow motives of products of
curves and Abelian varieties are finite dimensional \cite{K}, in
particular some non-ruled surfaces.  Moreover, if one
assumes Bloch-Beilinson conjecture, then all Chow motives
are finite dimensional, so Andr\'e  motivic zeta is conjecturally 
always rational.    

Let $X\subset \pp^n$ be a projective variety, and  $C_{p, d}(X)$ be
the Chow variety of $X$, which parametrizes the effective $p$-cycles
on $X$ with degree $d$.  In particular, the Chow variety  $C_{0,
  d}(X)$ of zero cycles of degree $p$  parametrizes a formal linear
combination of points $P_1+P_2+\cdots+P_d$ on $X$, hence the $d$-th
symmetric product $Sym^dX$ is canonically identified with
$C_{0,d}(X)$.  In \cite{E}, the first author proves that the formal
power series, the Euler-Chow series,  $\sum [\chi(C_{p, d}(X)]t^d$ is
rational for any simplicial projective toric variety, where $\chi(Y)$ is 
the Euler characteristic of $Y$.   

Observing these phenomena, it would be natural to ask if $\sum [C_{p,
  d}(X)]t^d$ is rational in $K(\ChM)[[t]]$.  In a previous paper by
the authors \cite{EK}, it was proved that if $n\geq 2$ then the series
\begin{equation}
\label{1}
\sum[C_{n-1,d}(\pp^n)]t^d=\sum_{d}\dfrac{1-[\AAA^1]^{\binom{n+d}{d}}}{1-[\AAA^1]}t^d  
\end{equation}
is irrational.  However, one may notice that in the irrational
power series above, if one takes the limit $[\AAA^1]\to 1$, then we
can recover the rationality of the Euler-Chow series.  So, we can
reformulate our last question as what happens if we change the
coefficients in the above series, in other words, what happens   if
one considers the above formal 
power series (\ref{1}) with coefficients the Chow motives of Chow varieties,
modulo the relation $[\AAA^1]=1$.  Is it rational? If it is, what
geometric or algebraic information can be read from it?  

The goal of this paper is to answer these questions as much as
possible. We are just going to get a glimpse of something that may 
morph to a deep subject, with a very hard and interesting questions.

Let us consider the $K$-ring $K(\ChM)$ of Chow motives, modulo $\AAA^1$
homotopy, namely we identify $[X\times \AAA^1]$ with $[X]$, we denote
this quotient by $K(\ChM)_{\AAA^1}$. Now,
define the Motivic Chow series by $MC_p(X):=\sum [C_{p, d}(X)]t^d$
with coefficients in $K(\ChM)_{\AAA^1}$. We will show this series is rational in many
cases. 

 (1) Picard Product formula, see Thm \ref{3-8}.  When ${\rm Pic}(X)\times {\rm
   Pic}(Y)={\rm Pic}(X\times Y)$, then we have
 $MC_{n-1}(X)\cdot MC_{m-1}(Y)=MC_{n+m-1}(X\times Y)$, where $n=\dim X$ and
 $m=\dim Y$.   In particular, when $X$ and $Y$ are curves, then
 $MC_0(X)$ and $MC_0(Y)$ are always rational by the result of
 Kapranov, hence for very general curves $X$ and $Y$, $MC_1(X\times
 Y)$ is rational with coefficients in $K(\ChM)_{\AAA^1}$.

(2) Torus action, see Thm \ref{4-3}: When the multiplicative group $\gm$ acts on $X$,
then, roughly speaking, $X$ is the disjoint union of the fixed locus $X^{\gm}$
and the free orbits $X_F$.  If the quotient $\overline{X}=X_F/\gm$ exists,
then the free orbit part looks like $\overline{X}\times \gm$, and its 
class in $K(\ChM)_{\AAA^1}$ is zero because $[\overline{X}\times
\gm]=[\overline{X}\times \AAA^1]-[\overline{X}]$.  
Thanks to Thomason's Torus generic slice theorem, see \cite[Prop. 4.10]{T},
one can justify this rough argument to show that $[X]=[X^{\gm}]$.
Moreover, it works for more general torus action to get $[X]=[X^T]$
when the algebraic torus $T$ acts on $X$.  This is a powerful tool to
compute the Motivic Chow series.  For example, one can show that the
Motivic Chow series of a toric variety is always rational.   

One important feature of the Motivic Chow series is that it can detect
very subtle geometric property of varieties.  For example, when one
blows up $\pp^2$ along $3$ points, then one can compute its cohomology
group, Chow group, higher Chow group, algebraic 
cobordims , or even Lawson cohomology group
without knowing the configuration of the $3$ points in $\pp^2$.  But
if one needs to 
calculate its Motivic Chow series, the answer depends on if the $3$
points are colinear or not, see examples \ref{gp} and \ref{n-points}.
 We should observe that example \ref{n-points} is new. It has not
 been able to be computed yet with the Euler-Chow series.  

We don't know if there exists any cohomology theory whose odd part and
even part correspond to the  numerator and the denominator of the
Motivic Chow series, but if such a cohomology theory exists, then it
will be a powerful tool to study algebraic varieties.  

\mbox{ }\newline
{\it Convention:}
Unless explicitly said,  throughout this paper we work in the category
of algebraic varieties over 
 a closed field of characteristic zero, denoted by $\kappa$. 
We also denote by $R$ any of the Grothendieck $K$-rings used in this
paper, for details see  
Definition \ref{1-9}. 

\mbox{ }\newline
{\it Acknowledgements:}  We are grateful for useful conversation and
personal communications with Jos\'e
de Jes\'us Malag\'on-L\'opez, Pedro Luis del Angel, Marc Levine and
James D. Lewis.  We also like to thank the referees for their
comments. 

\section{$K$-rings of categories}

\dfn{1-1} Let $\CCC$ be a category which has ``addition" (say
$X\disjoint Y$) and ``multiplication" (say $X\times Y$) such that  

\begin{enumerate}
\item[(1)] Addition is commutative and associative, namely $X\disjoint
  Y\simeq Y\disjoint X$, and $(X\disjoint Y)\disjoint Z\simeq
  X\disjoint (Y\disjoint Z)$. 
\item[(2)] Multiplication is commutative and associative, namely
  $X\times Y\simeq Y\times X$, and $(X\times Y)\times Z \simeq X\times
  (Y\times Z)$. 
\item[(3)] Distributive law holds, namely $X\times (Y\disjoint
  Z)\simeq (X\times Y)\disjoint (X\times Z)$. 
\end{enumerate}

An object $\unit\in \CCC$ is called {\em the multiplicative unit} if
$\unit\times X\simeq X$ for any $X\in \CCC$.  We always assume that the multiplicative unit exists.
Then we define {\em the $K$-ring of $\CCC$}, denoted as $K(\CCC)$, to
be the ring, generated by the objects of $\CCC$ (we write
the class of $X$ to be $[X]$), under the relations  

\begin{enumerate}
\item[(i)] if $X$ and $Y$ are isomorphic, then $[X]=[Y]$.
\item[(ii)] $[X\disjoint Y]-[X]-[Y]=0$.
\item[(iii)] $[X\times Y]-[X]\cdot [Y]=0$.
\item[(iv)] $[\unit]-1=0$.
\end{enumerate}

Sometimes, we consider the $K$-ring modulo extra relations, and we
denote the quotient as $K'(\CCC)$ when the extra relation is clear. 
\edfn

\rmk{1-2} When $\CCC={\mathcal FSets}$ is the category of finite sets,
usual disjoint union and multiplication satisfies the conditions, and
$K({\mathcal FSets})\simeq\zz$, by sending $[X]$ to $|X|$.   

 
 When $\CCC={\mathcal Top}$ is the category of ``good" topological
 spaces (e.g., finite simplicial complexes), and if we add the extra
 relation $[X]=[C]+[U]$ where $U\subset X$ is an open subset and
 $C\subset X$ the complement of $U$, then $K'({\mathcal Top})\simeq
 \zz$, by sending $[X]$ to $\chi(X)$, the Euler characteristic of $X$
 (see \cite{Y}). 
 
 We assume that a multiplicative unit $\unit$ exists, and its isomorphism class
 is unique.  If $\unit'$ is another multiplicative unit, then
 $\unit\simeq \unit\times \unit'\simeq \unit'$.   
 
 When $F: \CCC\to {\mathcal D}$ is a functor which preserves addition,
 multiplication, and the multiplicative unit if it exists, then $F$
 induces a ring homomorphism $K(\CCC)\to K({\mathcal D})$. 
 \ermk
 
 \dfn{1-3} Let $\CCC={\AVar}$ be the category of algebraic
 varieties over $\kappa$.  We define the addition in ${\AVar}$
 to be the disjoint union and the multiplication to be the product.
 We also consider the extra relation $[X]=[C]+[U]$ where $U\subset X$
 is an open subscheme and $C\subset X$ the complement of $U$ with the
 reduced scheme structure.  In this paper, $K'({\AVar})$
 always means $K({\AVar})$ modulo this extra relation. 
 
 Let ${\mathcal ChM}$ be the category of pure Chow motives.  Its
 objects are triple $(X, p, n)$ where $X$ is a smooth projective
 variety over $\kappa$, $p:X\vdash X$ an idempotent correspondence,
 and $n\in \Z$ (see  \cite{S}).  
  
  The category of Chow motives can be considered as
  ``the universal cohomology theory".
 
   In ${\mathcal ChM}$, we define
 addition as the direct sum and the multiplication as tensor product.   
 \edfn  
 
  \thm{1-4} (Bittner \cite{B}) $K'({\AVar})$ is canonically
  isomorphic to the ring  generated by smooth projective varieties,
  modulo the relation $[X]+[E]=[C]+[Bl_CX]$ where $C\subset X$ is a
  smooth closed subvariety, $Bl_CX\to X$ the blowing up along $C$, and
  $E\subset Bl_CX$ the exceptional divisor.  \ethm 
  
\cor{1-5} There is a ring homomorphism $K'(\AVar)\to K(\ChM)$ which
sends $[X]$ to $ch(X):=[(X, [\Delta_X], 0)]$, the Chow motive of $X$, where $\Delta_X\subset X\times X$ is the diagonal subvariety.  \ecor 

\begin{proof} It follows from the fact that in $\ChM$, we have
  $ch(X)\oplus ch(E)\simeq ch(C)\oplus ch(Bl_CX)$.
\end{proof} 

\rmk{1-6} Using this ring homomorphism $K'(\AVar)\to K(\ChM)$, we can
define $ch(X)$ to be the image of $[X]$, even when $X$ is neither
projective nor smooth.  For example, $ch(\AAA^1)=ch(\pp^1)-ch({\rm
  Pt})$, and when $X$ is a nodal rational curve, then
$ch(X)=ch(\AAA^1)$.  \ermk 

The following lemma will be useful in the next section.

\lem{1-7} Let $f: X\to Y$ be a proper morphism of reduced schemes,
which is bijective set theoretically.  Then $[X]=[Y]$  in $K'(\AVar)$.

In this case, as a morphism of topological spaces, $f$ is a homeomorphism.
\elem 

\begin{proof} First, we show that $f$ is a homeomorphism.  As $f_{\rm
    set}$ is continuous, we need to show that $f^{-1}_{\rm set}$ is
  also continuous.  As $f$ is proper, the image of a closed subset of
  $X$ by  $f_{\rm set}$  is closed, hence $f^{-1}_{\rm set}$ is
  continuous.  
  
  We may decompose $X$ and $Y$ into locally closed
  irreducible subschemes, so we may assume that $X$ and $Y$ are varieties.
  Let $\eta\in X$ be the generic point, then as $f$ is a homeomorphism,
  $f(\eta)\in Y$ is also the generic point.  The extension degree of
  $K(\eta)/K(f(\eta))$ must be $1$, otherwise $f$ cannot be bijective
  in characteristic $0$.  Then $f$ is birational, and there is an open
  subscheme $U\subset Y$ such that $f^{-1}(U)\to U$ is isomorphic, and
  in particular, $[f^{-1}(U)]=[U]$ in $K'(\AVar)$.  We need to show
  that $[X-f^{-1}(U)]=[Y-U]$, but $f_{|X-f^{-1}(U)}: X-f^{-1}(U)\to
  Y-U$ can be regarded as a proper morphism of reduced schemes, which
  is bijective set theoretically.  Hence by Noetherian induction, 
  we are reduced to the $0$-dimensional case, where Lemma
  is obvious.  \end{proof} 

\exm{1-8} When $X$ is a cuspidal rational curve, then its
normalization is a proper morphism which is bijective set
theoretically, and $[X]=[\pp^1]$ in $K'(\AVar)$.  \eexm 

\dfn{1-9} When $\CCC=\Hodge$ is the category of  the pure Hodge structures
with addition direct sum and multiplication tensor product, then
we can define $K(\Hodge)$.  For each pure Chow motive $M$, its cohomology
group $H^*(M)$ is well-defined and it has a Hodge structure, which
defines a ring homomorphism $K(\ChM)\to K(\Hodge)$.

When $\CCC=\Vect$ is the category of $\zz_2$-graded vector spaces,
with objects $V=V^{even}\oplus V^{odd}$.  We define the addition
as the direct sum and multiplication as the tensor product, with
$\Z_2$-grading following the usual convention.  Then we can
find $K(\Vect)\simeq \zz[\epsilon]/(\epsilon^2-1)$, where 
$[V^{even}\oplus V^{odd}]$ is sent to $\dim (V^{even})+\epsilon \dim (V^{odd})$.
There is a forgetful functor from $\Hodge$ to $\Vect$, which induces 
a natural ring homomorphism $K(\Hodge)\to K(\Vect)$.

From $K(\Vect)\simeq \Z[\epsilon]/(\epsilon^2-1)$, one can define a ring
homomorphism to $\Z$ by substituting $\epsilon=-1$.  The composition
$K'(\AVar)\to K(\ChM)\to K(\Hodge)\to K(\Vect)\to \Z$ sends $[X]$ to $\chi(X)$, 
the Euler characteristic.

For $\CCC=\AVar$, we can add one more relation, namely $[X\times \AAA^1]=[X]$,
which is called $\AAA^1$-homotopy relation.  We denote $K'(\AVar)$ modulo
$\AAA^1$-homotopy relation as $K'(\AVar)_{\AAA^1}$.  Similarly,
for $\CCC=\ChM$, we can also consider $\AAA^1$-homotopy relation, 
and we write $K(\ChM)$ modulo $\AAA^1$-homotopy relation as
$K(\ChM)_{\AAA^1}$. This is same as ignoring the dimension shifting,
namely $[(X, p, n)]\sim [(X, p, m)]$ for any $n, m\in \Z$.

 $\AAA^1$-homotopy relation for the Hodge structure
is equivalent to $[V(1)]=[V]$, ignoring the Tate twists.  We write 
$K(\Hodge)$ ignoring the Tate twists as $K(\Hodge)_{\hbox{\small Tate}}$.

Summarizing, we have the following commutative diagram for $K$-rings.

\xymatrix{ K'(\AVar) \ar[r] \ar@{->>}[d]
  \ar@{}[dr]|{\circlearrowright} &K(\ChM) \ar@{->>}[d]
  \ar@{}[dr]|{\circlearrowright} 
\ar[r] & K(\Hodge) \ar@{->>}[d] \ar[r]
\ar@{}[dr]|(.35){\circlearrowright} & K(\Vect) \ar[r] & \Z\\ 
K'(\AVar)_{\AAA^1} \ar[r] & K(\ChM)_{\AAA^1} \ar[r] &
K(\Hodge)_{\hbox{\small Tate}} \ar[ur] & \ } 
\edfn


\section{Chow varieties}

\dfn{2-1} Let $X\subset \pp^n$ be a projective variety, then we
denote by  $C_{p,d}(X)$ the Chow variety of $X$ which parametrizes
effective $p$-cycles of degree $d$ in $X$, see \cite{CvW}.  We
consider  $C_{p,d}(X)$ with its natural
reduced scheme structure.  We write $C_p(X):=\disjoint_{d\ge 0} C_{p, d}(X)$, and 
$B_p(X)$ the set of connected components of $C_p(X)$.  We also denote by $C_{p,\alpha} (X) $ the cycles of dimension $p$ in $X$ that are in the connected component $\alpha \in B_p(X).$  By the result of Hoyt (see the argument in Remark \ref{2-2}) $B_p(X)$ is independent  of the choice of the embedding $X\subset \pp^n$.\edfn 

\rmk{2-2} There is one-to-one correspondence between the closed points
of $C_{p, \alpha}(X)$ and the  effective cycles $\sum a_iV_i$ in the component of $\alpha\in B_p(X)$.   Hence the closed points of $C_p(X)$ corresponds 
one-to-one to effective $p$ cycles on $X$, and it is independent of
the choice of the embedding $X\subset \pp^n$ set theoretically, in
fact, as  a topological space with its Zariski topology, see \cite{H},
but Nagata showed that the scheme structure depends on the embedding
\cite{N}.  When $C_p(\xymatrix@1{X\,\, \ar@{^(->}^{\varphi}[r]& \pp^n})$
and $C_p(\xymatrix@1{X\,\, \ar@{^(->}^{\psi}[r]& \pp^m})$ are two Chow
varieties coming from different  embeddings into projective spaces,
then using the Chow variety $\widetilde C$ for the embedding
$\xymatrix@1{X\,\, \ar@{^(->}^(.35){\varphi\times \psi}[r] & \pp^n\times
  \pp^m\ \ar@{^(->}[r]& \pp^N}$ where $N=nm+n+m$ with
$\xymatrix@1{\pp^n\times \pp^m\ \ar@{^(->}[r]& \pp^N}$ the Segre
embedding, Hoyt proved that there are morphisms $\widetilde C$ to
$C_p(\xymatrix@1{X\ \ar@{^(->}^{\varphi}[r]& \pp^n})$ and
$C_p(\xymatrix@1{X\ \ar@{^(->}^{\psi}[r]& \pp^m})$ which are proper
bijective \cite{H}.  Hence by Lemma \ref{1-7}, we
have 
$$
[C_{p,\alpha}(\xymatrix@1{X\ \ar@{^(->}^{\varphi}[r]&
  \pp^n})]=[C_{p,\alpha}(\xymatrix@1{X\ \ar@{^(->}^{\psi}[r]& \pp^m})] \hbox{  in
} K'(\AVar)
$$ 
hence the class $[C_{p,\alpha}(X)]\in K'(\AVar)$ is independent
of the choice of the embedding into projective spaces. \ermk 

\rmk{2-3}
For a projective variety $X$, it is known that the Chow variety $C_{p, d}(X)$ is also
projective.  
The addition of effective cycles determines a proper morphism of schemes
$C_{p, d}(X)\times C_{p, e}(X)\to C_{p, d+e}(X)$, which gives an
abelian monoid object structure on $C_p(X)$, see Friedlander \cite[Prop. 1.8]{F}.  In
particular, its connected components $B_p(X)$ have an abelian monoid
structure.  The group completion of the monoid $B_p(X)$ is ${\rm
  CH}_p(X)/\stackrel{alg}{\sim}$, the Chow group modulo algebraic
equivalence \cite[Prop. 1.8]{F}. 
\ermk

\dfn{2-4} Let $X$ be a quasi-projective reduced scheme, and $\tilde
X\supset X$  be a projective completion, namely $\tilde X$ is a
projective scheme which has $X$ as an open subscheme.  Let $Y\subset
\tilde X$ be the complement of $X$ with reduced scheme structure, then
$C_p(Y)$ is a closed subscheme of $C_p(\tilde X)$.  The zero degree
part of $C_p(Y)$, namely $C_{p, 0}(Y)$ consists of one point
corresponding to $\phi$, and we define $C_p(Y)^+:=\disjoint_{d>0}C_{p,
  d}(Y)$, which is again a closed subscheme of $C_p(\tilde X)$.  We
define {\em the Chow variety of  the quasi-projective variety $X$ in
  $\tilde X$} to be the complement of the image of  $C_p(Y)^+\times
C_p(\tilde X)\to C_p(\tilde X)$ by the addition morphism in
$C_p(\tilde X)$, denoted as $C_p(X\subset \tilde X)$, which depends on
the embedding $\tilde X\in \pp^n$.  In other words, $C_p(X\subset \tilde X)$ 
is an open subscheme of $C_p(\tilde X)$ which consists of cycles 
$\sum a_i[V_i]$ with $V_i\cap X \not= \phi$. 

More generally, when $X$ is a locally closed subscheme of a projective
scheme $Z$, then take $\tilde X$ to be any closed subscheme of $Z$
which contains $X$.  Then $C_p(X\subset \tilde X)$ consists
of subset of $C_p(Z)$ which consists of cycles 
$\sum a_i[V_i]$ with $V_i\cap X\not=
\phi$, hence it is independent of the choice of $\tilde X$, 
so we define $C_p(X\subset Z)$ to be $C_p(X\subset \tilde X)$
for any such $\tilde X$.
\edfn 

\rmk{2-5}  By Remark \ref{2-2}, the class $[C_{p,\alpha}(X\subset Z)]\in
K'(\AVar)$ is independent of the embedding $\tilde X\subset \pp^n$.   

Set theoretically the closed points of $C_p(X\subset Z)$
correspond one-to-one to the effective $p$-cycles on $X$, and hence
it is independent of the choice of embedding $X\subset \pp^n$. 

But as a topological space, $C_p(X\subset Z)$ does depend on
the embedding.  For example, let $P\in \pp^2$ be a point, $X=\pp^2-P$,
 $Z_1=\pp^2$ and  $Z_2$  the blow-up of $Z_1$ along
$P$.  Then any two lines in $C_1(X\subset Z_1)$ are in the same
connected component, but in $C_1(X\subset Z_2)$, the lines
whose closure is through $P$ are in a different connected component as
the other lines. For the study of quasi-projective Chow varieties, interested readers can consult \cite{Li}.
\ermk


\lem{2-6} Let $X$ be a locally closed subscheme of a projective scheme $Y$, 
and $f: X \to Z$ a flat morphism of relative dimension $k$ to a quasi-projective scheme $Z$. 
 Let $\tilde Z$ be
any projective completion of $Z$, and $T\subset C_p(Z\subset \tilde Z)$ be any 
reduced locally closed subscheme.  Then there exists a stratification
$T=\disjoint T_i$ into locally closed  
subschemes of $T$ together with morphisms $\varphi_i:T_i\to
C_{p+k}(X\subset Y)$ such that for  
any point $t\in T_i$ corresponding  to a cycle $t=\sum n_iV_i$, the
image of the morphism $\varphi_i(t)$  
corresponds to $f^*(t)$, the flat pull-back by $f$ (see Fulton
\cite[Chap. 1]{Fu}).   

\elem

\begin{proof} By Noetherian induction, it is enough to show that there is a
non-empty open subscheme $U\subset T$ with a morphism $\varphi: U\to C_p(X\subset Y)$
such that for any $u\in U$, $\varphi(u)$ corresponds to $f^*(u)$ as a cycle.  

If ${\rm Im}(f)\cap T$ is not dense in $T$, then we can take $U\subset T$ so
that $U$ does not intersect with the image of $f$, and take the constant 
morphism $U\to \{\phi\}\in C_p(X\subset Y)$, so we may assume that
${\rm Im}(f)\cap T$ is  dense in $T$.  As flat morphism is open, by replacing $T$
with ${\rm Im}(f)\cap T$, we may assume that ${\rm Im}(f)\supset T$.

Let $U_1\subset T$ be an irreducible component of the smooth locus of
the reduced scheme $T$.   
By \cite[Theorem 1.4]{F},  the embedding $U_1\to C_p(Z\subset \tilde Z)$ 
corresponds to a cycle $\sum n_iW_i$ of $U_1\times \tilde Z$ such that 
for any point $u\in U$, the cycle $u\in C_p(Z\subset \tilde Z)$ 
corresponds to $\sum n_i[W_i\cap (\{u\}\times Z)]$.  In particular
the image of $W_i$ contains a dense open subscheme
of $U_1$.

We pull-back the cycle $\sum n_i[W_i\cap (U_1\times Z)]$ by the
flat morphism $1_{U_1}\times f$ to obtain $(1_{U_1}\times f)^*\sum n_i[W_i\cap (U_1\times Z)]=\sum m_j [\tilde W_j]$.
Let $\overline{W_j}$ be the closure of $\widetilde W_j$ in $U_1\times Y$.  
As the fiber of $\overline{W_j}\to U_1$ at the generic point is irreducible,
it is also irreducible in some open subscheme $U\subset U_1$ for all $j$, hence 
for $u\in U$, we have $\overline{W_j}\cap \{u\}\times
Y=\overline{\widetilde W_j\cap \{u\}\times Y}$.
Again by \cite[Theorem 1.4]{F},
there is a morphism $\varphi:U\to C_{p+k}(Y)$ such that for $u\in U$, the point
$\varphi(u)$ corresponds to the cycle $\sum m_j[\overline{W_j\cap \{u\}\times Y}]$, which
corresponds to $f^*(u)$, and in particular the image is in $C_{p+k}(X\subset Y)$.
We have constructed $\varphi: U\to C_{p+k}(X\subset Y)$.
\end{proof}

\dfn{2-7} Let $X$ be a locally closed subscheme of a projective variety $Y$,
and $f: X\to Z$ be a flat morphism of relative dimension $k$, and $Z\subset \tilde Z$
a projective completion.
For a subvariety $V\subset C_p(Z\subset \tilde Z)$ and $\alpha\in B_{p+k}(Y)$,
we define {\em the $\alpha$-component of $f^*[V]$} to be 
$$f^*[V]_{\alpha}:=\sum_i[\varphi_i(V_i)\cap C_{p+k, \alpha}(X\subset Y)]\in K'(\AVar)$$
where $V=\disjoint V_i$ is a stratification into locally closed
subschemes of $V$, and 
$\varphi_i: V_i\to C_{p+k}(X\subset Y)$  the morphism which maps $v\in V_i$
to the point corresponding to the cycle $f^*(v)$, as in Lemma \ref{2-6}.

For a linear combination of subvarieties, we extend $f^*$ linearly. 
\edfn

\lem{2-8} Let $X$ be a locally closed subscheme of a projective variety $Y,$ and $f: X\to Z$ be a surjective flat morphism of relative dimension $k$ to a quasi-projective variety $Z$, and $\alpha\in B_{p+k}Y$ a connected component of $C_{p+k}Y$.  Consider the class
$$\sum_{\beta\in B_p(\tilde Z)}f^*[C_{p, \beta}(Z\subset \tilde Z)_{\rm red}]_{\alpha}\in K'(\AVar)$$
for a projective completion $Z\subset \tilde Z$, then this class is independent of the choice of $\tilde Z$.\elem

\begin{proof} Let $Z\subset \tilde Z_i$, $i=1, 2$ be two projective completions, then 
both of $C_p(Z\subset \tilde Z_i)$ ($i=1, 2$) parametrize the effective $p$ cycles on $Z$, and hence there exists a common refinement $\tilde C_p\to C_p(Z\subset \tilde Z_i)$ $(i=1, 2)$ where $\tilde C_p$ is a disjoint union of locally closed subschemes of $C_p(Z\subset \tilde Z_i)$ ($i=1, 2$) with the reduced scheme structure so that the morphisms $\tilde C_p\to C_p(Z\subset \tilde Z_i)$ are bijective set theoretically.  Then $f^*[\tilde C_p]_{\alpha}$ makes sense for both of $C_p(Z\subset \tilde Z_i)$, well-defined independent of the choice of $\tilde Z_i$, and $f^*[\tilde C_p]=f^*[C_p(Z\subset \tilde Z_i)_{\rm red}]_{\alpha}$.\end{proof}

Lemma \ref{2-8} will be used in Definition \ref{3-11}, after we define  Motivic Chow series for quasi-projective varieties in the next section.

\section{Motivic Chow Series}

\dfn{3-1} An additive monoid $S$ is said to have {\em finite fiber}
when for any $s\in S$, the set $\{(a, b)\in S\times S|a+b=s\}$ is a
finite set.  We assume that additive monoid $S$ always have an
additive identity $0$. 
When $R$ is a commutative ring with a muliplicative unit $1$, and $S$
is an additive monoid with finite fiber, we define {\em the formal
  power series over $S$ with coefficient in $R$} to be the set of
functions from $S$ to $R$, written as $R[[S]]$, and we write an
element of $R[[S]]$ as $f=\sum_{s\in S}a_st^s$, where $f$ sends $s$ to
$a_s$.  Define the addition in $R[[S]]$ as the usual addition as
functions; $(\sum a_st^s)+(\sum b_st^s)=\sum (a_s+b_s)t^s$.  We define
the multiplication in $R[[S]]$ by convolution; $(\sum a_st^s)\cdot
(\sum b_st^s):=\sum_{s\in
  S}\left(\sum_{s_1+s_2=s}a_{s_1}b_{s_2}\right)t^s$, where the sum is finite,
  because $S$ has finite fiber.

A power series $f=\sum a_st^s\in R[[S]]$ is called {\em a polynomial}
when $a_s=0$ except for finitely many $s\in S$.  We denote the subring of all
polynomials as $R[S]$.   A polynomial $f=\sum a_st^s\in R[S]$ is
called {\em monic} when its coefficient of $t^0$ is $1$ (namely as a
function, $f$ sends $0$ to $1$).  A power series $f\in R[[S]]$ is
called {\em a rational function} if there exists a monic polynomial
$g\in R[S]$ such that $fg\in R[[S]]$ is a polynomial. 
\edfn

\rmk{3-2} When $X$ is a projective variety, the monoid of connected
components of the Chow variety $B_p(X)$ (see Definition \ref{2-1} and
Remark \ref{2-3}) has finite fiber, because for each fixed degree $d$,
the components of the Chow variety $C_p(X)$ corresponding to the
cycles with degree less than or equal to $d$ is a finite set.  When
$S_1\to S_2$ is a monoid homomorphism, then we have a ring
homomorphism $R[[S_1]]\to R[[S_2]]$.  For example, when $Y$ is a
closed subscheme of a projective variety $X$, then the monoid
homomorphism  $B_p(Y)\to B_p(X)$ induces $R[[B_p(Y)]]\to R[[B_p(X)]]$.
A ring homomorphism $R_1\to R_2$ canonically induces a ring
homomorphism $R_1[[S]]\to R_2[[S]]$.  \ermk

\dfn{3-3} When $X$ is a projective variety, then we define its {\em
  Motivic Chow Series of dimension $p$} to
be $$MC_p(X):=\sum_{\alpha\in B_p(X)} [C_{p, \alpha}(X)] t^{\alpha},$$ 
in $R[[B_p(X)]]$, 
where $R$ is any of the rings in Definition \ref{1-9}.

When $X$ is a quasi-projective variety with a fixed embedding into a
projective variety $X\subset \tilde X$, we define its {\em Motivic
  Chow Series of dimension $p$ in $\tilde X$} to be  
$$MC_p(X\subset \tilde X):=\sum_{\alpha\in B_p(\tilde X)}[C_{p,
  \alpha}(X\subset \tilde X)] t^{\alpha}.$$\edfn 

\exm{3-3-1} When $p=0$, the degree $d$ Chow variety $C_{0, d}(X)$ is
in the same class as $Sym^d(X)$ in $K'(\AVar)$, hence the Motivic Chow
series of $p=0$ is the Motivic zeta.  Kapranov \cite{Ka} proved that
when $X$ is a smooth projective curve, then $MC_0(X)$ is rational in
$K'(\AVar)[[B_0(X)]]\simeq K'(\AVar)[[t]]$, in the ring of formal
power series with one variable $t$.  When $X$ is a surface, Larsen
and Lunts  \cite{LL} proved that $MC_0(X)\in K'(\AVar)[[B_0(X)]] \simeq
K'(\AVar)[[t]]$ is rational if and only if $X$ is a ruled surface.  On
the other hand, using the notion of finite dimensionality of motives,
if the Chow motive of $X$ is finite dimensional, then $MC_0(X)\in
K(\ChM)[[B_0(X)]]\simeq K(\ChM)[[t]]$ is rational \cite{A}.  According
to Bloch-Beilinson Conjecture, all motives should be finite
dimensional, and when $X$ is a product of curves or Abelian variety,
their motives are finite dimensional, hence $MC_0(X)\in
K(\ChM)[[B_0(X)]]\simeq K(\ChM)[[t]]$ is rational \cite{K}.  For example, 
when $X=\pp^n$, we can write $[X]=\sum_{i=0}^n[{\mathbb A}^i] \in K'(\AVar)$ and $ch(\pp^n)=1+s+s^2+\cdots+s^n\in K(\ChM)$ with $s=ch(\AAA^1)=ch(\pp^1)-ch(Pt)$,  we have a rational expression
 $$MC_0(\pp^n)=\prod_{i=0}^n\dfrac{1}{1-s^it} \in K(\ChM)(B_0(\pp^n))=K(\ChM)(t).$$
But if we set $p>0$ for
$X=\pp^n$ with $n\ge 2$, we have $$MC_{n-1}(\pp^n)=\sum_{d=0}^{\infty}\frac{1-s^{\left(\begin{smallmatrix}d+n\\d\end{smallmatrix}\right)}}{1-s}t^d  \in K(\ChM)[[B_{n-1}]]=K(\ChM)[[t]],$$ 
which is
not rational by \cite{EK}.  

When $X$ is a smooth projective toric variety, then
$\displaystyle \sum_{\alpha\in H_2p(X, \Z)} \chi(C_{p,
  \alpha}(X))t^{\alpha}\in \zz[[H_{2p}(X, \Z)]]$ is rational by
\cite{E}.  Notice that there is a ring homomorphism
$K(\ChM)_{\AAA^1}\to \Z$ which sends $ch(X)$ to $\chi(X)$.
When $\{[V_1], [V_2], \ldots, [V_r]\}\subset {\rm CH}_*X$ is a $\zz$-basis of 
${\rm CH}_*X$, then one can find dense open subscheme $V_i\supset V_i^{\circ}\simeq \AAA^{d_i}$ for each $V_i$, so that $\displaystyle MC_0(X)=\prod_{i=1}^r \dfrac{1}{1-[V_i^{\circ}]t}$ in $K(\ChM)[[t]]$.

\eexm

\thm{3-4} (Localization) Let $p$ be a non-negative integer, $X$ a
projective variety, $U\subset X$ an open subvariety and $Y\subset X$
the complement of $U$ in $X$  with the reduced closed subscheme
structure.  Let $\overline{MC_p(Y)}$ be the image of $MC_p(Y)$ by the
ring homomorphism $R[[B_p(Y)]]\to R[[B_p(X)]]$  (see Remark
\ref{3-2}).  Then in $R[[B_p(X)]]$, we have 

$$\overline{MC_p(Y)}\cdot MC_p(U\subset X)=MC_p(X).$$\ethm

\begin{proof} Because there are ring homomorphisms from  $K'(\AVar)$ to 
all other $K$-rings, it is enough to show in the case
  $R=K'(\AVar)$. 

Let $\varphi: B_p(Y)\to B_p(X)$ be the natural morphism of monoids.
For each $\alphabar\in B_p(X)$, we define $\displaystyle C_{p,
  \alphabar}(Y):=\bigcup_{\varphi(\alpha)=\alphabar}C_{p, \alpha}(Y)$.
Set theoretically, each connected component $C_{p, \gamma}(X)$ of the
Chow variety of $X$ is a disjoint union of the images $C_{p,
  \alphabar}(Y)\times C_{p, \beta}(U\subset X)$ by the addition
morphism, with the index set $\{(\alphabar, \beta)\in B_p(X)\times
B_p(X)|\alphabar+\beta=\gamma\}$.  We denote the image of $C_{p,
  \alphabar}(Y)\times C_{p, \beta}(U\subset X)$ by $W_{\alphabar,
  \beta}\subset C_{p, \gamma}(X)$.  We will show that $W_{\alphabar,
  \beta}$ is a locally closed subscheme of $C_{p, \gamma}(X)$, and the
restriction of the addition morphism $C_{p, \alphabar}(Y)\times C_{p,
  \beta}(U\subset X)\to W_{\alphabar, \beta}$ is a proper bijection.
Because each irreducible cycle $V$ in $X$ is contained in one, and
only one of $C_p(Y)$ and $C_p(U\subset X)$, the way to decompose a
cycle in $C_p(X)$ into the sum of $C_p(Y)$ and $C_p(U\subset X)$ is
unique, from which the bijectivity follows. 

As $C_{p, \alphabar}(Y)\subset C_p(X)$ is closed, the image of $C_{p,
  \alphabar}\times C_{p, \beta}(X)$ in $C_{p, \gamma}(X)$ is closed.
As $C_{p, \beta}(U\subset X)\subset C_{p, \beta}(X)$ is open, the
image of $C_{p, \alphabar}(Y)\times \left(C_{p, \beta}(X)-C_{p,
    \beta}(U\subset X)\right)$ in $C_{p, \gamma}(X)$ by the addition
morphism is also closed.  By the following Lemma \ref{3-5}, if you
remove the image of\\ $C_{p, \alphabar}(Y)\times \left(C_{p,
    \beta}(X)-C_{p, \beta}(U\subset X)\right)$ in $C_{p, \gamma}(X)$
from the image of  $C_{p, \alphabar}\times C_{p, \beta}(X)$, then the
remaining subset is exactly $W_{\alphabar, \beta}$, hence
$W_{\alphabar, \beta}$ is locally closed. 

\lem{3-5} The inverse image of $W_{\alphabar, \beta}$ by the addition
morphism $C_{p, \alphabar}(Y)\times C_{p, \beta}(X)\to C_{p,
  \gamma}(X)$ is $C_{p, \alphabar}(Y)\times C_{p, \beta}(U\subset X)$.
\elem 

\begin{proof} (of Lemma \ref{3-5}) Each cycle corresponding to a point
  in $W_{\alphabar, \beta}$ can be written as $c+d$ with $c\in C_{p,
    \alphabar}(Y)$ and $d\in C_{p, \beta}(U\subset X)$.  Assume that
  it has another expression $c+d=c'+d'$ with $c'\in C_{p,
    \alphabar}(Y)$ and $d'\in C_{p, \beta}(X)$.  As $C_{p,
    \alphabar}(Y)\subset C_{p, \alphabar}(X)$ and $ C_{p,
    \alphabar}(X)$ is a connected component, we have $\deg c=\deg c'$.
  $c$ contains all the $Y$-supported cycles in $c+d=c'+d'$ with $c'$
  supported on $Y$, the cycle $c-c'$ is effective and degree $0$,
  hence $c=c'$, therefore $d=d'$, and $(c', d')\in C_{p,
    \alphabar}(Y)\times C_{p, \beta}(U\subset X)$.  (end of the proof
  of Lemma \ref{3-5}.)  \end{proof} 

From Lemma \ref{3-5}, it follows that the morphism $C_{p,
  \alphabar}(Y)\times C_{p, \beta}(U\subset X)\to W_{\alphabar,
  \beta}$ is the base extension of the proper morphism  $C_{p,
  \alphabar}(Y)\times C_{p, \beta}(X)\to C_{p, \gamma}(X)$ by the
inclusion $W_{\alphabar, \beta}\to C_{p, \gamma}(X)$, hence $C_{p,
  \alphabar}(Y)\times C_{p, \beta}(U\subset X)\to W_{\alphabar,
  \beta}$ is proper.   By Lemma \ref{1-7}, it follows that $[C_{p,
  \alphabar}(Y)]\times [C_{p, \beta}(U\subset X)]=[W_{\alphabar,
  \beta}]$. 

Set theoretically $C_{p, \gamma}(X)$ is a disjoint union of locally
closed subschemes $W_{\alphabar, \beta}$, in $K'(\AVar),$ hence in
$K(\ChM)$, and in $K(\ChM)_{\AAA^1}$. Therefore  $[C_{p,
  \gamma}(X)]=\sum_{\alphabar+\beta=\gamma} [C_{p,
  \alphabar}(X)]\times [C_{p, \beta}(U\subset X)]$.  Now we have 

\begin{multline*} 
\overline{MC_p(Y)}\times MC_p(U\subset
  X)=\\
= \left(\sum_{\alphabar\in B_p(X)}[C_{p,
      \alphabar}(Y)]t^{\alphabar}\right)\times \left(\sum_{\beta\in
      B_p(X)}[C_{p, \beta}(U\subset X)]t^{\beta}\right) \\ 
\hspace{-1.3cm}=  \sum_{\gamma\in B_p(X)}\left(\sum_{\alphabar+\beta=\gamma}[C_{p,
      \alphabar}(Y)]\times [C_{p, \beta}(U\subset X)]\right)
  t^{\gamma}\\\hspace{-6cm} = \sum_{\gamma\in B_p(X)}[C_{p, \gamma}(X)]
  t^{\gamma} \\ = MC_p(X).\mbox{\hspace{8.7cm}}
\end{multline*} 
\end{proof} 

\cor{3-6} Assume that $X$ is a projective variety, $Y\subset X$  is a
closed subscheme and $U\subset X$ its complement.  If two of following
series $\overline{MC_p(Y)},
MC_p(X)$ and $MC_p(U\subset X)$ are rational, then the other one is
also rational. \ecor 

\begin{proof} It follows from the fact that all $\overline{MC_p(Y)}, MC_p(X)$ and $MC_p(U\subset X)$ are  monic.    \end{proof}

\exm{3-7} When $Y=\{P_1, \ldots, P_r\}\subset X$ is a finite set of
$r$ elements in a connected scheme $X$, then
$\overline{MC_0(Y)}=\dfrac{1}{(1-t)^r}$, and $MC_0((X-Y)\subset X)\in
K(\ChM)[[t]]$ is rational if and only if the motive of $X$ is finite
dimensional.  In particular, when $U$ is any smooth curve with its
smooth completion $X$, $MC_0(U\subset X)\in K'(\AVar)[[t]]$ is
rational.  \eexm 

\rmk{3-15} It is interesting to see that when $Z=C\backslash \{P_1,
\ldots, P_{r+1}\}$ where $C$ is a smooth projective curve of genus $g$
and $r\ge 1$, then $Sym^n[Z]=0$ in $K'(\AVar)_{\AAA^1}$ if $n>2g+r-1$.
In particular, $MC_0(Z\subset C)$ is a polynomial. 

It is fun to see how it works concretely when $Z=\pp^1\backslash
\{a_1, a_2, \ldots, a_r,  \infty\}$.  We have isomorphism
$\AAA^r\simeq Sym^r(\AAA^1)$ by identifying $(b_1, b_2, \ldots,
b_r)\in\AAA^r$ with the $r$ roots of
$f(x)=x^r+b_1x^{r-1}+\cdots+b^{r-1}x+b_r$ in $Sym^r\AAA^1$.  Then
under this identification with the 
degree $r$ monic polynomials, $Sym^rZ$ corresponds to
$\{f(x)|f(a_1)\neq 0, f(a_2)\neq 0, \ldots, f(a_r)\neq 0\}$.  The
linear map $\AAA^r\to \AAA^r$ sending $(b_1, b_2, \ldots, b_r)$ to
$(f(a_1), f(a_2), \ldots, f(a_r))$ (where
$f(x)=x^r+b_1x^{r-1}+\cdots+b^{r-1}x+b_r$) is isomorphism by
Vandermonde, and by this new coordinate, we can see that
$[Sym^rZ]=[{\mathbb G}_m^r]=0$ in $K'(\AVar)_{\AAA^1}$. \ermk

\thm{3-8}  (Picard Product Formula) (1) Let $X$ and $Y$ be smooth
projective varieties such that the 
natural morphism gives the bijection $Pic(X)\times Pic(Y)\simeq
Pic(X\times Y)$.  Then we 
have 
$$
MC_{n-1}(X)\cdot MC_{m-1}(Y)=MC_{n+m-1}(X\times Y)\in
K'(\AVar)_{\AAA^1}[[B_{n+m-1}(X\times Y)]],
$$ 
where by the natural
morphism $Pic(X)\times Pic(Y)\simeq Pic(X\times Y)$, we identify
$B_{n-1}(X)\otimes B_{m-1}(Y)\simeq B_{n+m-1}(X\times Y)$, with
$n=\dim X$ and $m=\dim Y$. 

(2) The assumption $Pic(X)\times Pic(Y)\simeq Pic(X\times Y)$ holds if and
only if the only morphism of varieties $Pic^{\circ}(X)\to
Pic^{\circ}(Y)$ are the constant morphisms. 
 \ethm
 
 \begin{proof} (1)  Let $n=\dim X, m=\dim Y$ and ${\rm
     Pic}^{\circ}(X)$ and ${\rm Pic}^{\circ}(Y)$ 
 be the Picard varieties.  Let $NS(X)\subset H^2(X, \Z)$ be the Neron-Severi group,
 then the closed points of $NS(X)\times {\rm Pic}^{\circ}(X)$ corresponds one-to-one 
 to the element of ${\rm Pic}(X)$, the Pircard group of $X$, and
 similarly for ${\rm Pic}(Y)$.   
 
 For each $\alpha\in NS(X)$, the image of $C_{n-1}(X)$ in $\alpha\times {\rm Pic}^{\circ}(X)$
 is closed, and we denote it as $V_{\alpha}\subset \alpha\times {\rm Pic}^{\circ}(X)$.
 We have a stratification of $V_{\alpha}$ into locally closed subschemes $V_{\alpha}=\disjoint V_{\alpha, i}$
 such that the inverse image of $V_{\alpha, i}$ in $C_{n-1}(X)$ is $V_{\alpha, i}\times \pp^{d_i}$.  When
 $D$ is a divisor whose class is in $V_{\alpha, i}$, then we have $\dim H^0(X, \oo(D))=d_i+1$.
 Using these notations, we can write $MC_{n-1}(X)\in K(\ChM)_{\AAA^1}[[B_{n-1}(X)]]$ as
\begin{eqnarray*}MC_{n-1}(X)&=&\sum_{\alpha}\left(\sum_i[\pp^{d_i}]\cdot [V_{\alpha, i}]\right)t^{\alpha}\\
&=& \sum_{\alpha}\left(\sum_i(d_i+1)[V_{\alpha, i}]\right)t^{\alpha}. \end{eqnarray*}
We define $W_{\beta}=\disjoint W_{\beta, j}\subset \beta\times {\rm Pic}^{\circ}(Y)$ with the inverse
image of $W_{\beta, j}$ in $C_{m-1}(Y)$ be $W_{\beta, j}\times \pp^{e_j}$, we can write
$$MC_{m-1}(Y)=\sum_{\beta}\left(\sum_j(e_j+1)[W_{\beta, j}]\right) t^{\beta}.$$

When $D_1\in V_{\alpha, i}$ and $D_2\in W_{\beta, j}$ are effective
divisors, then $\pi_X^*D_1+\pi_Y^*D_2$ 
is in $V_{\alpha, i}\times W_{\beta, j}\subset NS(X)\times {\rm
  Pic}^{\circ}(X)\times NS(Y)\times {\rm Pic}^{\circ}(Y)$, 
and as $$H^0(X\times Y, \oo(\pi_X^*D_1+\pi_Y^*D_2))=H^0(X,
\oo(D_1))\otimes H^0(Y, \oo(D_2)),$$ 
its inverse image in $C_{n+m-1}(X\times Y)$ is $\pp^{(d_i+1)(e_j+1)-1}$, hence we have
\begin{multline*} 
MC_{n+m-1}(X\times Y)=\\
\hspace{-3.35cm}=\sum_{\alpha, \beta}\left(\sum_{i,
    j}[\pp^{d_ie_j+d_ie_j}]\cdot [V_{\alpha, i}]\cdot [W_{\beta,
    j}]t^{\alpha+\beta}\right)\\ 
\hspace{-2.35cm}=\sum_{\alpha, \beta}\left(\sum_{i, j}[(d_i+1)(e_j+1)]\cdot
  [V_{\alpha, i}]\cdot [W_{\beta, j}]t^{\alpha+\beta}\right)\\ 
\hspace{-.5cm}=\left(\sum_{\alpha}\left(\sum_i (d_i+1)[V_{\alpha, i}]\right)\right)\cdot 
\left(\sum_{\beta}\left(\sum_j (e_j+1)[W_{\beta, j}]\right)\right)\\
=MC_{2n-1}(X)\cdot MC_{2m-1}(Y)\hspace{6.3cm}
\end{multline*}

(2) First we assume that there is no non-constant morphisms ${\rm Pic}^{\circ}(X)\to {\rm Pic}^{\circ}(Y)$
as varieties. By restricting the line bundle on $X\times Y$ to $\{x\}\times Y$ and $X\times \{y\}$, one can easily recover
the preimage, so the natural morphism ${\rm Pic}(X)\times {\rm Pic}(Y)\to {\rm Pic}(X\times Y)$ 
is always injective, so we need to show that it is also surjective under the assumption.

 Take a line bundle $L$ on $X\times Y$.  Let $x\in X$ be a point,
and $L_Y:=L_{|\{x\}\times Y}$ be the restriction of $L$, considered as a line bundle on $Y$, and let
$M=(\pi_Y^*L_Y^{-1})\otimes L$, then $M_{|\{x\}\times Y}$ is the trivial line bundle.  
By the universality of ${\rm Pic}^{\circ}(Y)$, there exists a morphism $X\to {\rm Pic}^{\circ}(Y)$ 
such that $M\simeq (1_X\times \varphi)^*(\Poin_Y)\otimes \pi_X^*L_X$ for some line bundle $L_X$ on $X$,
where $\Poin_Y$ is the Poincar\'e line bundle on $Y\times {\rm Pic}^{\circ}(Y)$.  
The morphism $\varphi: X\to {\rm Pic}^{\circ}(Y)$ factors through $\tilde \varphi:{\rm Alb}(X)\to {\rm Pic}^{\circ}(Y)$,
and as ${\rm Pic}^{\circ}(X)$ is isogenous to ${\rm Alb}(X)$, there is a surjection
${\rm Pic}^{\circ}(X)\to {\rm Alb}(X)$.  So we have ${\rm Pic}^{\circ}(X)\to {\rm Alb}(X)\to {\rm Pic}^{\circ}(Y)$,
and by assumption, the composition ${\rm Pic}^{\circ}(X)\to {\rm
  Pic}^{\circ}(Y)$ is a constant morphism. 
Because ${\rm Pic}^{\circ}(X)\to {\rm Alb}(X)$ is surjective, the
morphism $\tilde \varphi$ is 
a constant morphism, and hence $\varphi: X\to {\rm Pic}^{\circ}(Y)$ is a constant morphism.
The image of $x\in X$ is $0$, because $M_{|\{x\}\times Y}$ is trivial, hence $\varphi$ is
a constant morphism to $0$, and $\varphi^{-1}\Poin_Y$ is the trivial line bundle.  
Therefore $M=\pi_X^*L_X$, and as $M=\pi_Y^*L_Y^{-1}\otimes L$, we conclude
$L=\pi_X^*L_X\otimes \pi_Y^*L_Y$.  We have shown that the natural morphism ${\rm Pic}(X)\times {\rm Pic}(Y)\to {\rm Pic}(X\times Y)$ is surjective, as required.

Conversely, assume that $\psi: {\rm Pic}^{\circ}(X)\to {\rm Pic}^{\circ}(Y)$ is
a non-constant morphism.   Take an isogeny ${\rm Alb}(X)\to {\rm Pic}^{\circ}(X)$
and let $\tilde\varphi: X\to {\rm Alb}(X)\to {\rm Pic}^{\circ}(X) \to {\rm Pic}^{\circ}(Y)$
be the composition, then it is a non-constant morphism, and $(1_X\times \tilde\psi)^*\Poin_Y$
cannot be written as $\pi_X^*L_X\times \pi_Y^*L_Y$, because for line bundle $\pi_X^*L_X\times \pi_Y^*L_Y$,
its restriction to $\{x\}\times Y$ is $L_Y$ for any $x\in X$, but the restriction of $(1_X\times \tilde\psi)^*\Poin_Y$
to $\{x\}\times Y$ corresponds to $\tilde\psi(x)$, which is not constant.  Hence we have
found a line bundle which is not in the image of ${\rm Pic}(X)\times {\rm Pic}(Y)$.
\end{proof}

\cor{3-9} Let $C_1, \ldots C_r$ be smooth projective curves such that
the morphisms of abelian varieties between their Jacobian varieties
$J(C_i)$'s are only zero morphisms.  Then $MC_{r-1}(C_1\times
C_2\times \cdots \times C_r)=\prod MC_0(C_i)$.  In particular, it is
rational in $K'(\AVar)_{\AAA^1}[[B_{r-1}(C_1\times\cdots\times C_r)]] $  \ecor 

\prodef{3-10} Let $X$ be a projective variety and $Y\subset X$ a
locally closed subscheme.  Let $\tilde Y\subset X$ be a closed
subscheme such that $Y\subset \tilde Y\subset X$, and $Y$ is an open
subscheme of $\tilde Y$.  We define $C_p(Y\subset X)$ to be
$C_p(Y\subset \tilde Y)$, considered as a locally closed subscheme of 
$C_p(X)$, and  $MC_p(Y\subset X)\in R[[B_p(X)]]$ to be
the canonical image of $MC_p(Y\subset \tilde Y)$ in $R[[B_p(X)]]$.
Then $C_p(Y\subset X)$ and $MC_p(Y\subset X)$ are independent of the
choice of $\tilde Y$.  When $X=\disjoint X_i$ is the set theoretic
decomposition of $X$ into locally closed subschemes, then we have
$\prod MC_p(X_i\subset X)=MC_p(X)$. \eprodef 

\begin{proof} $C_p(Y\subset\tilde Y)$ is an open subscheme of
  $C_p(\tilde Y)$, and as $C_p(\tilde Y)$ is a closed subscheme of
  $C_p(X)$, so we may regard $C_p(Y\subset \tilde Y)$ as the reduced
  locally closed subscheme of $C_p(X)$ whose points consist of linear
  combinations of irreducible varieties in $X$ whose scheme theoretic
  generic point is contained in $Y$.  We have found a description of
  $C_p(Y\subset X)$ without using $\tilde Y$, so it is independent of
  the choice of $\tilde Y$.  Let $\{\alpha\}$ be the set of connected
  components of $C_p(Y\subset X)$,  say $C_p(Y\subset X)=\disjoint
  C_{p, \alpha}(Y\subset X)$ as schemes.  Each $C_{p, \alpha}(Y\subset
  X)$ is contained in one connected component of $C_p(X)$, and for
  each $\beta\in B_p(X)$, we have $$MC_p(Y\subset X)=\sum_{\beta\in
    B_p(X)}\left(\sum_{C_{p, \alpha}(Y\subset X)\subset C_{p,
        \beta}(X)}[C_{p, \alpha}(Y\subset X)]\right) t^{\beta}$$ 
hence $MC_p(Y\subset X)$ is independent of the choice of $\tilde Y$.    

When $X=\disjoint X_i$ is a decomposition of $X$ into locally closed
subschemes, then at least one of $X_i$ is open, say $X_0$.  Then by
Theorem \ref{3-4}, letting $Y:=X-X_0$, we have
$\overline{MC_p(Y)}MC_p(X_0\subset X)=MC_p(X)$.  Now we can proceed by
induction on the number of locally closed subschemes, and the
observation that $\overline{MC_p(Y)}=MC_p(Y\subset X)$, together with
$MC_p(X)=\prod MC_p(X_i\subset X)$.  \end{proof}

\dfn{3-11} Let $X$ be a locally closed subscheme of a projective
variety $Y$, and $f: X\to Z$ a surjective flat  morphism of relative dimension
$k$, and $Z\subset \tilde Z$a projective completion.
Define $MC_p(f^*Z\to Y)$ to be
$$MC_p(f^*Z\to Y):=\sum_{\alpha\in B_p(Y)}\left(\sum_{\beta\in B_{p-k}(\tilde Z)}f^*[C_{p-k, \beta}(Z\subset \tilde Z)_{\rm red}]_{\alpha}\right)t^{\alpha}$$
See Definition \ref{2-7} for the notation $f^*[C_{p-k, \beta}(Z\subset \tilde Z)_{red}]_{\alpha}$.
This definition is independent of the choice of the projective completion $Z\subset \tilde Z$ by Lemma \ref{2-8}.
\edfn

\dfn{3-12} Let $X$ be a locally closed subscheme of a projective variety $Y$, $f:X\to Z$ a flat morphism of relative dimension $k$, and $Z\subset \tilde Z$ the projective completion.  
$$\xymatrix{X\ar@{^(->}[r] \ar[d]_{f} &Y\\ Z \ar@{^(->}[r] &\tilde Z}$$
This diagram is called {\em good configuration in dimension $p$} when $f$ is surjective and  there is a monoid homomorphism $\varphi: B_p(\tilde Z)\to B_{p+k}(Y)$ such that for any cycle $c\in C_{p, \alpha}(Z\subset \tilde Z)$, the flat pull-back $f^*(c)$ is in $C_{p+k, \varphi(\alpha)}(X\subset Y)$. \edfn

\pro{3-13} Let $X$ be a non-empty locally closed subscheme of a projective variety $Y$, $f:X\to Z$ a flat morphism of relative dimension $k$, and $Z\subset \tilde Z$ the projective completion.  When $\tilde Z$ is connected, then this configuration is a good configuration in dimension $0$ if and only if there exists $0\neq\beta\in B_p(Y)$ such that for any point $P\in Z$, we have $[f^{-1}(P)]\in C_{k, \beta}(X, Y)$.  

When our diagram is a good configuration in dimension $p$ with the
monoid homomorphism $\varphi: B_p(\tilde Z)\to B_{p+k}(Y)$, then the
following two conditions hold. 

\begin{enumerate}
\item[(1)] Let $\varphi_*: R[[B_p(\tilde Z)]]\to R[[B_{p+k}(Y)]]$ be
  the ring homomorphism induced by $\varphi$, with $R$ being $K'(\AVar)$,
  $K(\ChM)$ or $K(\ChM)_{\AAA^1}$, then we
  have $$\varphi_*(MC_p(Z\subset \tilde Z))=MC_{p+k}(f^*Z\to Y).$$ 
\item[(2)] When $MC_p(Z\subset \tilde Z)$ is rational, then
  $MC_{p+k}(f^*Z\to Y)$ is also rational. 
\end{enumerate}
\epro

\begin{proof} When $\tilde Z$ is connected, we have $B_0(\tilde
  Z)\simeq \Z_{\ge 0}$ with $C_{0, d}\tilde Z\simeq Sym^d\tilde Z$, as
  topological spaces.  Each point $P\in \tilde Z$ can be regarded as
  an element of $\tilde Z\simeq C_{0, 1}(\tilde Z)$, and
  $f^*[P]:=[f^{-1}(P)]$ being in the same component $C_{k, \beta}(Y)$
  is a necessary condition for our diagram to be a good configuration.
  Conversely, assume that $[f^{-1}(P)]\in C_{0, \beta}(X\subset Y)$
  for all point $P\in Z$.  As $X$ is non-empty, at least for one $P\in
  Z$, $f^{-1}(P)$ is not an empty set, so $\beta$ cannot be the
  connected component of the empty set.  In particular, $f$ is
  surjective.  For any cycle $P_1+\cdots+P_d\in C_{0, d}(Z\subset
  \tilde Z)$, we have $f^*(P_1+\cdots +P_d)=\sum_{i=1}^df^*(P_i)\in
  C_{k, d\beta}(X\subset \tilde Y)$, so we can define $\varphi:
  B_0(\tilde Z)\to B_k(Y)$ to be $\varphi(d):=d\beta$ to make our
  diagram a good configuration in dimension $0$. 

When our diagram is a good configuration in dimension $p$, then by
definition of $MC_{p+k}(f^*Z\to Y)$, we have 
$$MC_{p+k}(f^*Z\to Y) = \sum_{\beta\in B_p(Y)}\left(\sum_{\alpha\in
    B_{p}(\tilde Z)}f^*[C_{p, \alpha}(Z\subset \tilde
  Z)]_{\beta}\right)t^{\beta}$$ 
where $$f^*[C_{p, \alpha}(Z\subset \tilde
Z)]_{\beta}=\sum_i[\varphi_i(V_i)\cap C_{p+k, \beta}(X\subset Y)]$$
with $C_p(X\to Y)=\disjoint V_i$ is a stratification into locally
closed subschemes, $\varphi_i: V_i\to C_p(X\to Y)$ sends $c\in
C_p(Z\subset \tilde Z)$ to $f^*(c)\in C_p(X\subset Y)$.  Because our
diagram is a good configuration, for each $\alpha\in B_p(\tilde Z)$,
any $c\in C_{p, \alpha}(Z\subset \tilde Z)$ is sent to
$C_{p,\varphi(\alpha)}(X\subset Y)$ injectively (because $f$
is surjective), hence $f^*[C_p(Z\subset \tilde
Z)]_{\beta}=\sum_{\varphi(\alpha)=\beta}[C_p(Z\subset \tilde Z)]$.
Therefore we have 
\begin{eqnarray*}MC_{p+k}(f^*Z\to Y)&=&\sum_{\beta\in
    B_p(Y)}\left(\sum_{\varphi(\alpha)=\beta}[C_{p,\alpha}(Z\subset
    \tilde Z)]\right)t^{\beta}\\ &=&\sum_{\alpha\in B_p(\tilde
    Z)}[C_{p,\alpha}(Z\subset \tilde
  Z)]t^{\varphi(\alpha)}\end{eqnarray*} 
which is exactly the image of $\sum_{\alpha\in B_p(\tilde Z)}[C_{p,
  \alpha}(Z\subset \tilde Z)]t^{\alpha}=MC_p(Z)$ by $\varphi_*$.  

When for a monic polynomial $g\in R[B_p(\tilde Z)]$, $MC_p(Z)\cdot g$ is a polynomial, then $\varphi_*(MC_p(Z)\cdot g)=\varphi_*(MC_p(Z))\cdot \varphi_*(g)=MC_{p+k}(f^*Z\to Y)\cdot \varphi_*(g)$, and as $\varphi_*(g)$ is a monic polynomial in $R[B_p(Y)]$, we have shown that $MC_{p+k}(f^*Z\to Y)$ is rational.\end{proof}

\cor{3-14} Let $X\subset Y$ be a locally closed subscheme of a projective variety $Y$, and $f: X\to Z$ a flat morphism with relative dimension $k$ to $Z\hookrightarrow C$, an open curve with the projective completion $C$.  When $C\backslash Z=\{P_1, P_2, \ldots, P_r\}$ and there exists $\beta\in B_k(Y)$ such that  $f^*[P]\in C_{k, \beta}(X\subset Y)$ holds for any $P\in Z$, then 
$MC_p(f^*Z\to Y)$ is rational for any $p$.  When $r\ge 2$, then $MC_k(f^*Z\to Y)\in K'(\AVar)_{\AAA^1}$ is a polynomial.  When moreover $Z$ is a rational curve, then $MC_k(f^*Z\to Y)=(1-t^{\beta})^{r-2}\in K'(\AVar)_{\AAA^1}$.  \ecor

\begin{proof} By Proposition \ref{3-13}, our configuration is a good
  configuration in dimension $0$.  Also it is easy to verify that it
  is a good configuration in dimension $1$, and for all dimensions
  (except for $0$ and $1$, the condition is empty).  By Kapranov's
  theorem \cite{Ka}, $MC_0(C)$ is rational, with denominator
  $(1-t)^2\in K'(\AVar)_{\AAA^1}$ where $t$ is the class of a point.
  By Localization Thoerem \ref{3-4}, we have $MC_0(Z\subset
  C)\overline{MC_0(\{P_1, \ldots, P_r\})}=MC_0(C).$  As $\overline{MC_0(\{P_1, \ldots,
  P_r\})}=\frac{1}{(1-t)^r}$, we see that $MC_0(Z\subset C)$ is
  rational, and is a polynomial when $r\ge 2$.   
When $C$ is $\pp^1$, we have $MC_0(Z\subset C)=(1-t)^{r-2}$ and by
Proposition \ref{3-13} (2), $MC_k(f^*Z\to Y)=(1-t^{\beta})^{r-2}$.
For the other $p$, the equalities $MC_1(Z\subset
C)=\frac{1}{1-t^{[C]}}$ and $MC_p(Z\subset C)=1$, $p\not= 0,1$  are easy to
check.  \end{proof}

\section{Torus action}

\lem{4-2} Let $T={\mathbb G}_m^n$ be a torus.  Assume $T$ acts on a reduced scheme
$X$ with the fixed point locus $X^T$, then we have $[X]=[X^T]$ in $K'(\AVar)_{\AAA^1}$.  \elem 

\begin{proof} By Thomason's Torus generic slice theorem
  \cite[Prop. 4.10]{T}, $X$ can be decomposed into locally closed
  subschemes $X=\disjoint X_i$ such that the stabilizer on $X_i$ is
  constant, say $\gm^{d_i}\simeq H_i\subset T$, the quotient 
  $X_i/(T/H_i)=\overline{X_i}$ exists and $X_i\simeq
  \overline{X_i}\times (T/H_i)$.  When $d_i<n$, then $[X_i]=0$
  in $K(\AVar)_{\AAA^1}$.  In fact, we can write $X_i\simeq Y_i\times \gm$
  where $Y_i=\overline{X_i}\times {\mathbb G}_m^{n-d_i+1}$, then 
  $[X_i]=[Y_i]\times [\gm]=[Y_i]\times [\AAA^1]-[Y_i]=0$.
Hence we have
  $[X]=\sum_{H_i=T}[X_i]=[X^T]$.  \end{proof} 

\thm{4-3} Assume that $X$ is a reduced scheme on which the torus
$T={\mathbb G}_m^n$ acts.  By Thomason's Torus generic slice theorem
\cite[Prop. 4.10]{T}, we can write $X=\disjoint X_i$, a decomposition
into locally closed subschemes $X_i$ so that the stabilizer on $X_i$
is $H_i\simeq {\mathbb G}_m^{d_i}\subset T$, $X_i/(T/H_i)=:\overline{X_i}$
exists and $X_i\simeq \overline{X_i}\times (T/H_i)$.  Let $\pi_i:X_i\to \overline{X_i}$ be
the quotient map, then we have

$$MC_p(X)=\prod_i MC_p(\pi_i^*\overline{X_i}\to X).$$
\ethm

\begin{proof} We have
\begin{eqnarray*} MC_p(X)&=& \sum_{\alpha\in B_p(X)}[MC_{p, \alpha}(X)]t^{\alpha}\\
&=& \sum_{\alpha\in B_p(X)}[MC_{p, \alpha}(X)^T]t^{\alpha}\end{eqnarray*}
by Lemma \ref{4-2}.  A cycle $\sum n_iV_i$ is a fixed point by $T$ 
if and only if each $V_i$ is fixed by $T$.  An irreducible variety
$V_i$ is fixed by $T$ if and only if as a cycle $V_i\in
MC_p(\pi_j^*\overline{X_j}\to X)$ for  
some $j$.  

So for a point $\sum n_iV_i\in C_{p, \alpha}(C)$, it is in $C_{p, \alpha}^T$ if and only
if  it is in the locus $\prod_kf^*[\overline{X}_{i_k}\subset \tilde
  X_{i_k}]_{\beta_{j_k}}$ with $\sum \beta_{j_k}=\alpha$, for some
projective completion $\overline{X}_{i_k}\subset \tilde X_{i_k}$ for
each $k$, 
and once we fix the projective completion $\overline{X}_{i_k}\subset \tilde X_{i_k}$, 
 this decomposition is unique, so $[C_{p, \alpha}(X)]$ is the
 coefficient of $t^{\alpha}$ in 
$\prod_i MC_p(\pi_i^*\overline{X}_i\to X)$.\end{proof}

\dfn{4-4} For an irreducible subvariety $V\subset X$, we define {\em the degree of $V$},
denoted as $\deg V$ to be the class of $V$ in $B_p(X)$.  \edfn

\cor{4-5} 
When $X$ is a toric variety, then we have
$$MC_p(X)=\prod_{V} \frac{1}{(1-t^{\deg V})}$$
where $V$ runs over all the $p$-dimensional orbits.\ecor

\begin{proof} As a toric variety $X$ has finitely many orbits by the torus action, 
so in order to use the formula of Theorem \ref{4-3}, 
we may take $\{X_i\}$ to be the
set of orbits, and so all the quotients $\overline{X_i}=P_i$ are single points. 
Then as $$MC_p(\pi_i^*P_i\to X)=\left\{\begin{array}{cc} 1& \hbox{$p\not=\dim X_i$}\\
\frac{1}{1-t^{\deg X_i}} & \hbox{$p=\dim X_i$}\end{array}\right.,$$
by Theorem \ref{4-3} we have
$$MC_p(X)=\prod_{\dim X_i=p}\frac{1}{1-t^{\deg X_i}}$$ \end{proof}

\rmk{referee} When $p=0$, the formula $MC_p(X)=\prod_{V} \frac{1}{(1-t^{\deg V})}$
is valid in $K(\ChM)[[B_0(X)]]=K(\ChM)[[t]]$ without ${\mathbb A}^1$-homotopy relation.
\ermk

\exm{gp} Let $P_i\in \pp^2$, $(i=1,2,3)$ be three points on $\pp^2$, and we assume that
they are not colinear.  Let $X$ be the blow-up of these three points.
We can assume that  these points are $(1:0:0), (0:1:0), (0:0:1)$, and
then by the usual 
torus action, one can regard $X$ as a complete toric variety, whose
fan has $1$-dimensional 
cones spanned by $\{(1,0), (1,1), (0,1), (-1,0), (-1, -1), (0, -1)\}$ respectively.  The irreducible divisors corresponding
to $(1, 1), (-1, 0), (0,-1)$ are the exceptional divisors, and we denote them as
$E_1, E_2, E_3$, with $P_i\in \pp^2$ the image of $E_i$ for each $i$.
Let $\overline{L_i}\subset \pp^2$ be the line through $P_j, P_k$ where $i, j, k$ are chosen as 
$\{i, j, k\}=\{1, 2, 3\}$, and $L_i\subset X$ be the strict transform of $\overline{L_i}$, then
$L_1, L_2, L_3$ are the irreducible $T$-divisors corresponding to
$(-1, -1), (1, 0), (0,1)$.  The effective divisors on $X$ are exactly the
linear combination of $[L_1], [L_2], [L_3], [E_1], [E_2], [E_3]$, with 
non-negative coefficients.  The linear equivalence relation 
is generated by $[L_i]+[E_j]=[L_j]+[E_i]$.  We write $t_i=t^{[L_i]}$
and $s_i=t^{[E_i]}$, then by Corollary \ref{4-5}, the motivic Chow series of 
$X$ is $$MC_1(X)=\frac{1}{(1-t_1)(1-t_2)(1-t_3)(1-s_1)(1-s_2)(1-s_3)}.$$
These $6$ variables have the relation $t_is_j=t_js_i$.  

Let $[L_0]:=[L_i]-[E_i]$ (which is the same class for $i=1, 2, 3$), and let $t_0:=t^{[L_0]}$,
then $L_0$ is not effective, but $\Z[t_1, t_2, t_3, s_1, s_2, s_3]\subset \Z[t_0, s_1, s_2, s_3]$
with $t_i=t_0s_i$, and the $4$ variables $\{t_0, s_1, s_2, s_3\}$ are algebraically independent, 
and we have another expression of the motivic Chow series
as $$MC_1(X)=\frac{1}{(1-t_0s_1)(1-t_0s_2)(1-t_0s_3)(1-s_1)(1-s_2)(1-s_3)}.$$
\eexm

\exm{n-points}
 In this example, we compute the motivic Chow series of the blow-up
of $\pp^2$ along colinear $r$ points.

Let $\overline{L}:=\{(0:y:z)\}\subset \pp^2$ be the $y$-axis, $\{P_1,
P_2, \ldots, 
P_r\}\subset \overline{L}$ be $r$ distinct points on $\overline{L}$.
Let $\pi:X\to \pp^2$ be the blow-up 
along $\{P_1, \ldots, P_r\}$, and $L\subset X$ be the strict transform
of $\overline{L}$, with $\tilde P_i:=\pi^{-1}(P_i)\cap L$. 
Define the action of $\gm$ on $\pp^2$ by $\lambda (x: y: z):=(\lambda
x:y:z)$, then this action extends to 
an action on $X$ uniquely.  In this action, the fixed locus $X^{\gm}$
consists of $L$ and finitely many 
points, more precisely, $Q_0:=\pi^{-1}(1:0:0)$ and $Q_i\in
E_i\backslash L$ where $E_i$ is the exceptional divisor  
over $P_i$. Consider the morphism $\displaystyle \lim_{\lambda\to 0}
\lambda P: X\to X^{\gm}$, then 
$\psi^{-1}(L-\{\tilde P_1, \ldots, \tilde P_r\})\cap(X\backslash
X^{\gm})\simeq (L-\{\tilde P_1, \ldots, \tilde P_r\})\times \gm$. 
Define $X_0:=\psi^{-1}(L-\{\tilde P_1, \ldots, \tilde P_r\})\cap(X\backslash X^{\gm})$.  
The other free orbits are $X_i:=\psi^{-1}Q_i\cap (X\backslash X^{\gm})$ and $E_i^{\circ}:=E_i\cap (X\backslash X^{\gm})$.
Write $t^{[L]}=:t_0$ and $t^{[E_i]}=:s_i$.  Also let $p_i: X_i\to X_i/\gm$ and $q_j: E_j^{\circ}\to E_j^{\circ}/\gm$ be the quotient
morphisms
Then we can write the motivic Chow series of $X$ as
\begin{eqnarray*}MC_1(X)&=&MC_1(X^{\gm})\prod_{i=0}^r MC_1(p_i^*X_i/\gm\to X)\prod_{j=0}^r MC_0(q_j^*E_i^{\circ}/\gm\to X)\\ &=& \frac{1}{1-t_0}\cdot (1-t_0s_1s_2\cdots s_r)^{r-2}\prod _{i=1}^r\frac{1}{1-t_0s_1s_2\cdots \hat{s_i}\cdots s_r}\prod_{j=1}^r \frac{1}{1-s_i}\end{eqnarray*}

where we used Corollary \ref{3-14} for the computation of $MC_1(p_0^*X_0/\gm\to X)$.  In particular when $r=3$, comparing with Example \ref{gp}, we notice that the motivic Chow series depend on the configuration of the center of the blow-up.

In the case of $p=0$ in contrast, it is known that the Motivic zeta $MC_0(X)=\sum [Sym^n ch(X)]t^n$ depends only on the Chow motive of $X$.  One can blow-up $\pp^2$ along a point $r$ times, whose Chow motive is isomorphic to $ch(X)$, and whose motivic zeta is $$MC_0(X)=\dfrac{1}{(1-t)(1-[X^{\circ}]t)(1-[\tilde L^{\circ}]t)}\prod_{i=1}^r \dfrac{1}{1-[E_i^{\circ}]t}$$ where $\tilde L\subset X$ is the strict transform of a line $L\subset \pp^2$ which does not intersect with the blow-up centers, $\tilde L^{\circ}=\tilde L\backslash \{P\}$ for a $\kappa$-valued point $P\subset \tilde L$, and each $E_i$ is the strict transform of the exceptional divisor of the $i$-th blow-up,    $E_i^{\circ}=E_i\backslash \{P_i\}$ with $P_i\in E_i$ any $\kappa$-valued point, and $[X^{\circ}]=[X]- ([\tilde L]+\sum [E_i^{\circ}])$. \eexm

\rmk{final} \cite{KKT} says that the motivic Chow series may not be rational for $\pp^2$ blown-up along too many general points.  Can we expect the rationality of the motivic Chow series at least for varieties defined over $\ff_1$, which might be a natural assumption, because we are taking the limit $\AAA^1\to 1$ (counting the number of the points over $\ff_1$)? \ermk


\vskip 5mm

E. Javier Elizondo\\ \ \ \ \ Instituto de Matem\'aticas, Universidad Nacional Aut\'onoma de M\'exico.\\ \ \ \ \ 04510 M\'exico DF, M\'exico \\ \ \ \ \ e-mail:
javier@math.unam.mx

\vskip 5mm

Shun-ichi Kimura\\ \ \ \ \ Department of Mathematics, Hiroshima
University\\  \ \ \ \  Higashi-Hiroshima, Kagamiyama 1-3-1, 739-8526,
JAPAN\\  \ \ \ \ e-mail: kimura@math.sci.hiroshima-u.ac.jp


\begin{thebibliography}{k} 
\bibitem{A} {\sc Y. Andr\'e}. {\it Motifs de dimension finite
    (d'apr\`es S. I. Kimura, P. Sullivan...)}, in S\'eminaire Bourbaki
  Vol. 2003/2004 Ast\'erisque, vol. 299, Soci\'et\'e Math\'ematique de
  Frace, 2005, Exp. 929, 115--145 
\bibitem{B} {\sc F. Bittner}. {\it The universal Euler characteristic
    for varieties of characteristic zero}. Compos. Math. 140 1011 --
  1032 (2004) 
\bibitem{CvW} {\sc W-L Chow, B. L. van der Wareden}. {\it Zur
    algebraischen Geometrie IX. \"Uver zugeordnete Formen und
    algebraische Systeme von algebraischen
    Mannigfaltigkeiten}. Math. Ann. 113 (1937) 692 -- 704 
\bibitem{D} {\sc B. Dwork}. {\it On the rationality of the zeta function of an algebraic variety},
American J. Math. (1960) 631--648
\bibitem{E} {\sc E. J. Elizondo}. {\it The Euler series of restricted
    Chow varieties}, Compos. Math. 94 (1994) 297--310 
\bibitem{EK} {\sc E. J. Elizondo, S. Kimura}. {\it Irrationality of
    Motivic series of Chow varieties}, Math. Z. 263 (2009) 27--32 
\bibitem{F} {\sc E. Friedlander}. {\it Algebraic cycles, Chow
    varieties, and Lawson homology}, Compositio Math. 77 (1991) 55-93 
\bibitem{Fu} {\sc W. Fulton}. {\it Intersection Theory}, Springer (1984)
\bibitem{H} {\sc W. Hoyt}. {\it On the Chow bunches for different
    projective embeddings of a complete variety}, Amer. J. Math. 88
  (1966) 273--278 
\bibitem{Ka}
{\sc M. Kapranov.}
{\it The elliptic curve in the S-duality theory
and Eisenstein series for Kac-Moody group}s, 
preprint(arXiv: math.AG/0001005).
\bibitem{K} {\sc S. Kimura}. {\it Chow groups are finite dimensional,
    in some sense}, Mathematische Annalen {\bf 331}, 173-201 (2005) 
 \bibitem{KKT} {\sc S. Kimura, S. Kuroda and N. Takahashi} {\it Irrationality of Certain Euler-Chow sereis}, in preparation
\bibitem{LL}
{\sc M. Larsen, V.A. Lunts.}
{\it Rationality criteria for motivic zeta functions,} 
Compos. Math. {\bf 140}(2004), no. 6, 1537--1560.
\bibitem{Li}{\sc P. Lima-filho.}{\it Lawson Homology for
    quasiprojective varieties},   Compositio Mathematica, tome 84,
  No. 1, p. 1-23.
\bibitem{Mu} {\sc D. Mumford.} {\it Abelian Varieties} 2nd ed. Tata
  Institute of Fundamental Research, Oxford University Press (1974) 
\bibitem{N} {\sc M. Nagata.} {On the normality of the Chow variety of
    positive $0$-cycles of degree $m$ in an algebraic variety},
  Memoirs of the College of Science, Kyoto, Series A. vol. XXIX (1955)
  165--176 
\bibitem{S} {\sc T. Scholl.} {\it Classical Motives}, in Proceedings
  of Symposia in Pure Mathematics Volume 51 (1994), Part 1, 163--187 
\bibitem{T} {\sc R. W. Thomason.} {\it Comparison of equivariant
    algebraic and topological $K$-theory}, Duke Math 53 (1986)
  795--825 
\bibitem{W} {\sc A. Weil}. {\it Numbers of solutions of equations in finite fields},
Bulletin of AMS 55 (1949) 497--508
\bibitem{Y} {\sc S. Yokura.} {\it Motivic Characteristic classes}, in Topology of stratified spaces,  Cambridge Univ. Press (2011) 375--418
\end{thebibliography}
\end{document}